\documentclass[11pt]{amsart}
\usepackage{amsmath}
\usepackage{amscd}
\usepackage{graphics}
\usepackage{latexsym}
\oddsidemargin .in
\theoremstyle{plain}
\newtheorem{Thm}{Theorem}
\newtheorem{Lem}[Thm]{Lemma}
\newtheorem{Cor}[Thm]{Corollary}
\newtheorem{Prop}[Thm]{Proposition}
\newtheorem{Rm}{Remark}
\newtheorem{Ex}{Example}
\newtheorem{Def}{Definition}

\theoremstyle{remark}

\def\C{{\mathbb C}}

\def\R{{\mathbb R}}

\def\P{{\mathbb P}}

\def\Z{{\mathbb Z}}

\def\O{{\mathbb O}}
\def\H{{\mathbb H}}

\def\CL{\mathcal L}

\def\s2x{\hbox{$S^2 \times S^2$}}

    \def\sqr#1#2{{\vcenter{\hrule height.#2pt
            \hbox{\vrule width.#2pt height#1pt \kern#1pt
            \vrule width.#2pt}\hrule height.#2pt}}}
    \def\square{\mathchoice\sqr67\sqr67\sqr{2.1}6\sqr{1.5}6}

\def\qed{~\hfill$\square$}

\begin{document}

\title[]{Mirror Duality via $\mathbf{G_2}$ and $\mathbf{Spin(7)}$ Manifolds}
\author{Selman Akbulut and Sema Salur}
\thanks{First named author is partially supported by NSF grant DMS 0505638}
\keywords{mirror duality, calibration}
\address{Department  of Mathematics, Michigan State University, East Lansing, MI, 48824}
\email{akbulut@math.msu.edu }
\address {Department of Mathematics, University of Rochester, Rochester, NY, 14627 }
\email{salur@math.rochester.edu } \subjclass{53C38,  53C29, 57R57}
\date{\today}

\begin{abstract}
The main purpose of this paper is to give a construction of  certain ``mirror dual'' Calabi-Yau submanifolds inside of a $G_2$ manifold.  More specifically, we explain how to assign a $G_2$
manifold $(M,\varphi, \Lambda)$, with the calibration 3-form
$\varphi$ and an oriented  $2$-plane field $\Lambda$, a pair of
parametrized tangent bundle valued 2 and 3-forms of $M$. These
forms can then be used  to define different complex and
symplectic structures on certain  6-dimensional subbundles of
$T(M)$. When these bundles are integrated they give mirror CY
manifolds. In a similar way, one can define mirror dual $G_2$
manifolds inside of a $Spin(7)$ manifold $(N^8, \Psi)$. In case
$N^8$ admits an oriented  $3$-plane field, by iterating this
process we obtain Calabi-Yau submanifold pairs in $N$ whose
complex and symplectic structures determine each other via the
calibration form of the ambient $G_2$ (or $Spin(7)$) manifold.

\end{abstract}

\maketitle

\setcounter{section}{0}
\vspace{-0.3in}

\section{Introduction}

Let $(M^7,\varphi)$ be a $G_2$ manifold with the calibration 3-form
$\varphi$. If $\varphi$ restricts to be the volume form of an oriented
3-dimensional submanifold $Y^3$, then $Y$ is called an associative
submanifold of $M$. Associative submanifolds are very interesting
objects as they behave very similarly to holomorphic curves of
Calabi-Yau manifolds.

\vspace{.05in}

In \cite{as}, we studied the deformations of associative
submanifolds of $(M, \varphi)$ in order to construct Gromov-Witten
like invariants. One of our main observations was that oriented
$2$-plane fields on $M$ always exist by a theorem of Thomas
\cite{t}, and by using them one can split the tangent bundle
$T(M)={\bf E}\oplus {\bf V}$  as an orthogonal direct sum of an
associative $3$-plane bundle ${\bf E}$ and a complex $4$-plane
bundle $ {\bf V}$. This allows us to define ``complex associative
submanifolds'' of $M$, whose deformation equations may be  reduced
to the Seiberg-Witten equations, and hence we can assign local
invariants to them, and assign various invariants to $(M,\varphi,
\Lambda)$, where $\Lambda$ is an oriented $2$-plane field on $M$.
It turns out that these Seiberg-Witten equations on the
submanifolds are restrictions of global equations on $M$.

\vspace{.05in}

In this paper, we explain how the geometric structures on $G_2$
manifolds with oriented $2$-plane fields $(M,\varphi, \Lambda )$
provide complex and symplectic structures  to certain
6-dimensional subbundles of $T(M)$.  When these bundles integrated we obtain a pair of  Calabi-Yau manifolds whose complex and symplectic structures are remarkably related to each other.
We also study examples of Calabi-Yau manifolds which
fit nicely in our mirror set-up. Later, we do similar
constructions for $Spin(7)$ manifolds with oriented $3$-plane
fields. We then explain how these structures lead to the
definition of ``dual $G_2$ manifolds'' in a $Spin(7)$ manifold, with their own dual Calabi-Yau submanifolds.

\vspace{.07in}

{\em Acknowledgements:} We would like to thank R.Bryant and
S.Gukov for their valuable comments.

\vspace{.1in}

\section{Associative and Complex distributions of a $G_2$ manifold}

Let us  go through quickly over the basic definitions about $G_2$
manifolds.  The main references are the two foundational papers
\cite{hl} and \cite{b1}, as well as \cite{s}, \cite{b2},
\cite{bs},  and \cite{j}. We also need some properties introduced
in \cite{as}. Now let $\O=\H\oplus l \H= \R^8$  be the octonions
which is an $8$ dimensional division algebra generated by  $<1, i,
j,  k, l, li ,lj, lk> $, and let $im \O =\R^7$ be the imaginary
octonions with the cross product operation $\times:  \R^7\times
\R^7 \to \R^7$, defined  by $u\times v=im(\bar{v}.u)$. The
exceptional Lie group $G_{2}$ is  the linear automorphisms of  $im
\O$ preserving this cross product operation, it can also be
defined in terms of the orthogonal $3$-frames in $\R^7$:

\begin{equation*}
G_{2} =\{ (u_{1},u_{2},u_{3})\in (im \O)^3\;|\: \langle u_{i},u_{j} \rangle=\delta_{ij},  \;
\langle u_{1} \times u_{2},u_{3} \rangle =0 \;  \}.
\end{equation*}

Another very useful definition popularized  in \cite{b1} is the subgroup of $GL(7,\R)$ which fixes a particular $3$-form
$\varphi_{0} \in \Omega^{3}(\R ^{7})$.  Denote
$e^{ijk}=dx^{i}\wedge dx^{j} \wedge dx^{k}\in \Omega^{3}(\R^7)$,
then

$$ G_{2}=\{ A \in GL(7,\R) \; | \; A^{*} \varphi_{0} =\varphi_{0}\; \}. $$
\vspace{-0.15in}
\begin{equation}
\varphi _{0}
=e^{123}+e^{145}+e^{167}+e^{246}-e^{257}-e^{347}-e^{356}.
\end{equation}

\vspace{.05in}

{\Def A smooth $7$-manifold $M^7$ has a {\it $G_{2}$ structure} if
its tangent frame bundle reduces to a $G_{2}$ bundle.
Equivalently, $M^7$ has a {\it $G_{2}$ structure} if there is  a
3-form $\varphi \in \Omega^{3}(M)$  such that  at each $x\in  M$
the pair $ (T_{x}(M), \varphi (x) )$ is  isomorphic to $(T_{0}(
\R^{7}), \varphi_{0})$ (pointwise condition). We call $(M,\varphi)$ a manifold with $G_2$ structure.}

\vspace{.05in}

A  $G_{2}$ structure $\varphi$ on $M^7$ gives an orientation
$\mu \in \Omega^{7}(M)$ on $M$,   and $\mu$  determines  a
metric $g=g_{\varphi }= \langle \;,\;\rangle$ on $M$, and a cross product
structure $\times$  on the tangent bundle of $M$ as follows: Let
$i_{v}=v\lrcorner $ be the interior product with a vector $v$, then

\begin{equation}
\langle u,v \rangle=[ i_{u}(\varphi ) \wedge i_{v}(\varphi )\wedge \varphi  ]/6\mu .
\end{equation}
\begin{equation}
\varphi (u,v,w) = \langle u\times v,w \rangle .
\end{equation}

\vspace{.05in}

{\Def A manifold with $G_{2}$ structure $(M,\varphi)$  is called a
{\it $G_{2}$ manifold} if   the holonomy group of the Levi-Civita
connection (of the metric $g_{\varphi }$) lies inside of  $G_2 $.
Equivalently  $(M,\varphi)$ is a $G_{2}$ manifold if $\varphi $ is
parallel with respect to the metric $g_{\varphi }$, that is
$\nabla_{g_{\varphi }}(\varphi)=0$; which is equivalent to
$d\varphi=0 $, $\;d(*_{g_{\varphi}}\varphi)=0$. Also equivalently,  at each point
$x_{0}\in M$ there is a chart  $(U,x_{0}) \to (\R^{7},0)$ on which
$\varphi $ equals to $\varphi_{0}$ up to second order term, i.e.
on the image of $U$
$\varphi (x)=\varphi_{0} + O(|x|^2)$}.

\vspace{.05in}

{\Rm One important class of $G_2$ manifolds are the ones obtained
from Calabi-Yau manifolds. Let $(X,\omega, \Omega)$ be a complex
3-dimensional Calabi-Yau manifold with K\"{a}hler form $\omega$
and a nowhere vanishing holomorphic 3-form $\Omega$, then
$X^6\times S^1$ has holonomy group $SU(3)\subset G_2$, hence is a
$G_2$ manifold. In this case $ \varphi$= Re $\Omega + \omega
\wedge dt$. Similarly, $X^6\times \mathbb{R}$ gives a noncompact
$G_2$ manifold.}

\vspace{.05in}

{\Def Let $(M, \varphi )$ be a $G_2$ manifold. A 4-dimensional
submanifold $X\subset M$ is called {\em coassociative } if
$\varphi|_X=0$. A 3-dimensional submanifold $Y\subset M$ is called
{\em associative} if $\varphi|_Y\equiv vol(Y)$; this condition is
equivalent to the condition $\chi|_Y\equiv 0$,  where $\chi \in \Omega^{3}(M,
TM)$ is the  tangent bundle valued 3-form defined by the
identity:}
\begin{equation}
\langle \chi (u,v,w) , z \rangle=*\varphi  (u,v,w,z)
\end{equation}
The equivalence of these  conditions follows from  the `associator equality' of  \cite{hl}
\begin{equation*}
\varphi  (u,v,w)^2 + |\chi (u,v,w)|^2/4= |u\wedge v\wedge w|^2
\end{equation*}

Similar to the definition of $\chi$ one can define a tangent
bundle 2-form, which is just the cross product of $M$
(nevertheless viewing it as a $2$-form has its advantages).

\vspace{.05in}

{\Def Let $(M, \varphi )$ be a $G_2$ manifold. Then $\psi \in
\Omega^{2}(M, TM)$ is the tangent bundle valued 2-form defined by
the identity:}
\begin{equation}
\langle \psi (u,v) , w \rangle=\varphi  (u,v,w)=\langle u\times v , w
\rangle
\end{equation}

\vspace{.05in}

Now we have two useful properties  from \cite{as}, the first
property basically follows from definitions, the second property
fortunately applies when the first property fails to give anything
useful.

\vspace{.05in}

{\Lem (\cite{as})  To any $3$-dimensional submanifold $Y^3\subset
(M,\varphi)$,  $\chi$ assigns a normal vector field, which
vanishes when $Y$ is associative.}

{\Lem(\cite{as}) To any associative manifold $Y^3\subset (M,\varphi)$ with a non-vanishing oriented $2$-plane field,  $\chi$ defines a complex structure on  its normal bundle (notice in particular that any coassociative submanifold $ X\subset M$ has an almost complex structure if its normal bundle has a non-vanishing section).}

\begin{proof}
Let $L\subset \R^7$ be an associative $3$-plane,
that is $\varphi_{0}|_{L}=vol(L)$. Then for every pair of
orthonormal vectors $\{u,v\}\subset L$, the form $\chi$ defines a
complex structure on the orthogonal $4$-plane $L^{\perp}$, as
follows:   Define $j: L^{\perp} \to L^{\perp}$ by
\begin{equation}
j(X)=\chi(u,v,X)
\end{equation}
This is well defined i.e. $j(X)\in L^{\perp}$, because when $ w\in L$ we have:
$$ \langle \chi(u,v,X),w \rangle=*\varphi_{0}(u,v,X,w)=-*\varphi_{0}(u,v,w,X)=\langle \chi(u,v,w),X \rangle=0$$

\noindent Also $j^{2}(X)=j(\chi(u,v,X))=\chi(u,v,\chi(u,v,X))=-X$. We can check the last equality by taking an orthonormal basis $\{ X_{j}\}\subset L^{\perp}$ and calculating
\begin{eqnarray*}
\langle \chi(u,v,\chi(u,v,X_{i})),X_{j}\rangle &=&*\varphi_{0}(u,v,\chi(u,v,X_{i}),X_{j})=-*\varphi_{0}(u,v,X_{j},\chi(u,v,X_{i}))\\
&=&- \langle \chi(u,v,X_{j}),\chi(u,v,X_{i})\rangle =-\delta_{ij}
\end{eqnarray*}

The last equality holds since the map $j$ is orthogonal, and the
orthogonality  can be seen by polarizing the associator equality, and by noticing $\varphi_{0}(u,v,X_i)=0$. Observe that the
map $j$ only depends on the oriented $2$-plane $\Lambda=<u,v>$ generated by $\{u,v\}$ (i.e. it only depends on the complex structure on
$\Lambda $). \end{proof}

{\Rm Notice that Lemma 1 gives an interesting flow on the $3$-dimensional submanifolds of $G_2$ manifolds
$f: Y \hookrightarrow (M,\varphi )$ (call $\chi$-flow), described by:
$$\frac{\partial}{\partial t} f=\chi(f_{*}vol \;(Y) )$$
  For example, by \cite{bs} the total space of the spinor bundle $Q^7\to S^3$ (with $\C^2$ fibers) is a $G_2$ manifold, and the zero section $S^3\subset Q$ is an associative submanifold. We can imbed any homotopy $3$-sphere $\Sigma^{3} $  into $Q$ (homotopic to the zero-section). We conjecture that $\chi $- flow on $\Sigma\subset Q$, takes $\Sigma$ diffeomorphically onto the zero section. Note that, since any $S^3$ smoothly unknots in $S^7$ it is not possible to produce  quick counterexamples  by tying local knots; and an affirmative answer  gives $\Sigma \cong S^3$. }

  \begin{figure}[ht]  \begin{center}
\includegraphics{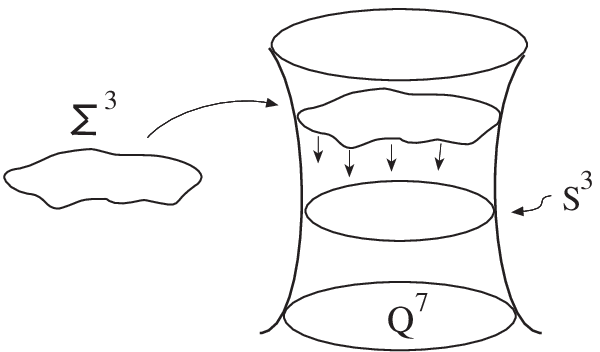}   \caption{}
\end{center}
   \end{figure}

Finally, we need some identities from \cite{b2} (also see \cite{k}) for $(M^7,\varphi)$, which follow from local calculations by using the definition (1). For
$\beta \in \Omega^{1}(M)$ we have:
\begin{equation}
|\;\varphi\wedge \beta \;|^2=4|\beta|^2, \; \;\mbox{and}\; \; |*\varphi\wedge \beta \;|^2=3|\beta|^2 ,  \\
 \end{equation}
  \begin{equation}
  (\xi \lrcorner \; \varphi ) \wedge \varphi =2*(\xi \lrcorner \;\varphi )\;,\;\mbox{and}\; \;
*\;[ \; *( \beta \wedge *\varphi )\wedge * \varphi \;]=3\beta ,
 \end{equation}
 \begin{eqnarray}
\beta \times (\beta \times u)=-|\beta |^2u+\langle \beta,u\rangle \beta ,
\end{eqnarray}
 where $*$ is the star operator. Let $\xi$ be a vector field on any Riemannian manifold $(M, g )$, and $\xi^{\#}\in \Omega^{1}(M)$ be its dual $1$-form, i.e. $\xi^{\#}(v)=\langle \xi,v\rangle$.  Then for
$\alpha \in \Omega^{k}(M)$:
\begin{eqnarray}
*(\xi \lrcorner \; \alpha )&=&(-1)^{k+1}(\xi^{\#} \wedge *\alpha) .
\end{eqnarray}

\vspace{.05in}

\section{Mirror duality in $G_2$ manifolds}

\vspace{.05in}

On a local chart of a $G_2$ manifold $(M, \varphi )$, the form $\varphi $ coincides with the form $\varphi_{0} \in \Omega^{3} (\R^7)$  up to quadratic terms, we can express the corresponding tangent valued forms $\chi$ and $\psi$ in terms of $\varphi_{0}$ in local coordinates. More generally, if $e_1,...e_7$ is any local orthonormal frame and $e^1,..., e^7$ is the dual frame, from definitions we get:

 \begin{equation*}
\begin{aligned}
\chi=&\;\;(e^{256}+e^{247}+e^{346}-e^{357})e_1\\
&+(-e^{156}-e^{147}-e^{345}-e^{367})e_2\\
&+(e^{157}-e^{146}+e^{245}+e^{267})e_3\\
&+(e^{127}+e^{136}-e^{235}-e^{567})e_4\\
&+(e^{126}-e^{137}+e^{234}+e^{467})e_5\\
&+(-e^{125}-e^{134}-e^{237}-e^{457})e_6\\
&+(-e^{124}+e^{135}+e^{236}+e^{456})e_7.\\
\end{aligned}
\end{equation*}

\vspace{.1in}

\begin{equation*}
\begin{aligned}
\psi=&\;\;(e^{23}+e^{45}+e^{67})e_1\\
&+(e^{46}-e^{57}-e^{13})e_2\\
&+(e^{12}-e^{47}-e^{56})e_3\\
&+(e^{37}-e^{15}-e^{26})e_4\\
&+(e^{14}+e^{27}+e^{36})e_5\\
&+(e^{24}-e^{17}-e^{35})e_6\\
&+(e^{16}-e^{25}-e^{34})e_7.\\
\end{aligned}
\end{equation*}

\vspace{.1in}

The forms $\chi$ and $\psi$ induce complex and symplectic
structures on certain subbundles of $T(M)$ as follows:
Let $\xi$ be a  nonvanishing vector field of $M$. We can define a symplectic $\omega_{\xi}$ and a complex  structure  $J_{\xi}$ on the $6$-plane bundle $V_{\xi}:=\xi^{\perp}$ by
\vspace{.05in}
\begin{equation}
\omega_\xi=\langle \psi, \xi \rangle \;\;\;\mbox{and} \;\; J_{\xi}(X)=X\times \xi.
\end{equation}

Now we can define
\vspace{.05in}
\begin{equation}
 \textup{Re}\; \Omega_{\xi} = \varphi|_{V_{\xi}}
\;\;\mbox{and}\;\;\; \textup{Im}\; \Omega_{\xi} = \langle \chi, \xi \rangle.
\end{equation}

\vspace{.1in}

\noindent In particular $\omega_{\xi}=\xi \lrcorner \; \varphi$, and $\textup{Im}\;\Omega_{\xi} =\xi \lrcorner \; *\varphi $. Call $\Omega_{\xi} =  \textup{Re}\; \Omega_{\xi} +i \;\textup{Im}\; \Omega_{\xi}$. The reason for defining these is to pin down a Calabi-Yau like structure on any $G_2$ manifold. In case $(M, \varphi)=CY \times S^1$ these quantities are related to the ones in Remark 1. Notice that when $\xi \in {\bf E}$ then  $J_{\xi}$ is an extension of $J$ of Lemma 2 from the $4$-dimensional bundle ${\bf V}$ to the $6$-dimensional bundle $V_{\xi}$.

\vspace{.1in}

 By choosing different directions, i.e. different $\xi$, one can
find the corresponding complex and symplectic structures. In
particular we will get two different complex structures if we
choose $\xi$ in the associative subbundle ${\bf E}$ (where $\varphi$ restricts to be 1), or if we choose $\xi$  in the complementary subbundle ${\bf V}$, which we will call the coassociative subbundle. Note that $\varphi$ restricts to zero on the coassociative subbundle.

\vspace{.1in}

In local coordinates, it is a straightforward calculation that by choosing $\xi=e_i$ for any $i$, from equations (11) and (12), we can easily obtain the corresponding structures $\omega_{\xi}$, $J_{\xi}$,
$\Omega_{\xi}$.
For example, let us assume that $\{e_1,e_2,e_3\}$ is the local orthonormal basis for the associative bundle ${\bf E}$,  and
$\{e_4,e_5,e_6,e_7\}$ is the local orthonormal basis for the coassociative bundle ${\bf V}$. Then if we choose $\xi=e_3=e_1\times e_2$
then we get $\omega_{\xi}= e^{12}-e^{47}-e^{56}$
and Im $\Omega_{\xi}=
e^{157}-e^{146}+e^{245}+e^{267}$. On the other hand, if we choose
$\xi=e_7$ then $\omega_{\xi}= e^{16}- e^{25}-e^{34}$
and Im $\Omega_{\xi}=
-e^{124}+e^{135}+e^{236}+e^{456}$ which will
give various symplectic and complex structures on the bundle $V_{\xi}$.

\vspace{.05in}

\subsection {A useful example}

$\;$
\vspace{.1in}

Let us take a Calabi-Yau
6-torus $\mathbb {T}^6=\mathbb{T}^3\times  \mathbb{T}^3$, where $\{e_1,e_2,e_3\}$ is the basis for one $\mathbb{T}^3$ and $\{e_4,e_5,e_6\} $ the basis for the other (terms expressed with a slight abuse of notation). We can take the product $M=\mathbb {T}^6\times S^1$ as the corresponding $G_2$ manifold with the calibration 3-form
$\varphi=e^{123}+e^{145}+e^{167}+e^{246}-e^{257}-e^{347}-e^{356}
$, and with the decomposition $T(M)={\bf E}\oplus {\bf V}$, where ${\bf E}=\{e_1,e_2,e_3\}$ and ${\bf V}=\{e_4,e_5,e_6,e_7\}$. Now, if we choose $\xi=e_{7}$, then $V_{\xi}=<e_1,...,e_6>$ and the symplectic form is $\omega_{\xi}=
e^{16}-e^{25}-e^{34}$,  and the complex structure is
 \[
J_{\xi} =\left(
\begin{array}{ccc}
{\bf e_1}  & \mapsto  &-e_6   \\
{\bf e_2} & \mapsto  & e_5  \\
 {\bf e_3} & \mapsto & e_4
\end{array}
\right)
\]

\vspace{.05in}

\noindent and the complex valued
$(3,0)$ form is $\Omega _{\xi }=(e^1+ie^6)\wedge
(e^2-ie^5)\wedge(e^3-ie^4)$;  note that this is just
$\Omega _{\xi }=(e^1-iJ_{\xi}(e^1))\wedge
(e^2-iJ_{\xi}(e^2))\wedge(e^3-iJ_{\xi}(e^3)) $.

\vspace{.05in}

On the other hand, if we choose $\xi '=e_{3}$ then
 $V_{\xi '}=< e_1,..,\hat{e}_{3},..,e_7>$ and the symplectic form  is
 $\omega_{\xi '}= e^{12}-e^{47}-e^{56}$ and the complex structure is
 \[
J_{\xi '} =\left(
\begin{array}{ccc}
{\bf e_1 } & \mapsto  &- {\bf e_2 }  \\
 e_4 & \mapsto  & e_7  \\
  e_5 & \mapsto & e_6
\end{array}
\right)
\]
\vspace{.1in}

\noindent Also
 $\Omega_{\xi '} =(e^1+ie^2)\wedge
(e^4-ie^7)\wedge(e^5-ie^6)$, as above this can be expressed  more tidily as  $\Omega_{\xi '}=(e^1-iJ_{\xi'}(e^1))\wedge
(e^4-iJ_{\xi'}(e^4))\wedge(e^5-iJ_{\xi'}(e^5))$.
In the expressions of $J$'s  the basis of associative bundle ${\bf E}$ is indicated by bold face letters to indicate the differing complex structures on $\mathbb{T}^6$. To sum up:  If we choose $\xi$ from the coassociative bundle ${\bf V}$ we get the complex structure which decomposes the 6-torus as $\mathbb{T}^3\times \mathbb{T}^3$.  On the other hand if we choose $\xi$ from the associative bundle ${\bf E}$ then the induced complex structure on the $6$-torus corresponds to the decomposition as $\mathbb{T}^2\times \mathbb{T}^4$. This is the phenomenon known as ``mirror duality''. Here these two $SU(3)$ and $SU(2)$ structures are different but they come from the same $\varphi$ hence they are dual. These examples suggests the following definition of ``mirror duality'' in $G_2$ manifolds:

\vspace{.05in}

{\Def Two Calabi-Yau manifolds are mirror pairs of each other, if
their complex structures are  induced from the same calibration
3-form in a $G_2$  manifold. Furthermore we call them strong
mirror pairs if their normal vector fields $\xi$ and $\xi'$ are
homotopic to each other through nonvanishing vector fields. }

\vspace{.05in}

{\Rm In the above example of $CY\times S^1$, where $CY=
\mathbb{T}^6$, the calibration form $\varphi=$Re $\Omega
+\omega\wedge dt $ gives Lagrangian tori fibration in $X_{\xi} $
and complex tori fibration in $X_{\xi'}$. They are different
manifestations of  $\varphi $ residing on one higher dimensional
$G_2$  manifold $M^7$. In the next section this correspondence
will be made precise.}

\vspace{.1in}

In Section 4.2 we will discuss a more general notion of mirror
Calabi-Yau manifold pairs, when they sit in different $G_2$
manifolds, which are themselves mirror duals of each other  in a
$Spin(7)$ manifold.

\subsection{General setting}

$\;$
\vspace{.1in}

Let $(M^7, \varphi , \Lambda)$ be a manifold with a $G_2$ structure and a non-vanishing  oriented  $2$-plane field. As suggested in \cite{as}  we can view $(M^7, \varphi)$ as an analog of a symplectic manifold, and the $2$-plane field $\Lambda $ as an analog of a complex structure taming $\varphi$. This is because $\Lambda $ along with
$\varphi $ gives the associative/complex bundle splitting $T(M)={\bf E}_{\varphi,\Lambda}\oplus {\bf V}_{\varphi,\Lambda}$.  Now, the next object is a choice of a non-vanishing unit vector field $\xi \in \Omega^{0}(M,TM)$, which gives a codimension one distribution $V_{\xi}:= \xi^{\perp}$ on $M$, which is  equipped with the structures $(V_{\xi},  \omega_{\xi}, \Omega_{\xi}, J_{\xi})$ as given by (11) and (12).

\vspace{.05in}

Let  ${\xi}^{\#}$ be the dual $1$-form of $\xi$. Let  $e_{\xi^{\#}}$  and $i_{\xi}=\xi \lrcorner $ denote the exterior and interior product operations on differential forms. Clearly
$e_{\xi^{\#}}\circ i_{\xi }+i_{\xi }\circ e_{\xi^{\#}} =id  $.
\begin{equation}
\varphi =e_{\xi^{\#}}\circ i_{\xi }(\varphi )+i_{\xi }\circ e_{\xi^{\#}}(\varphi )=\omega_{\xi}\wedge  \xi^{\#}+Re \; \Omega_{\xi}.
\end{equation}
This is just  the decomposing of the form $\varphi $ with respect to $\xi \oplus \xi^{\perp}$.
Recall that the  condition that the distribution $V_{\xi}$ be integrable (the involutive condition which implies $\xi^{\perp}$ comes from a foliation)  is given by:
\begin{equation}
 d{\xi}^{\#}\wedge {\xi}^{\#}=0.
 \end{equation}
 Even when  $V_{\xi}$ is not integrable, by \cite{th} it is homotopic to a foliation. Assume $X_{\xi }$ be a page of this foliation; for simplicity assume this 6-dimensional manifold is smooth.

\vspace{.05in}

It is clear from definitions that  $J_{\xi}$ is an almost complex
structure on $X_{\xi}$. Also the $2$-form $\omega_{\xi }$ is
non-degenerate on $X_{\xi}$, because from (2) we can write
\begin{equation}
\omega_{\xi }^3=(\xi \lrcorner \;\varphi )^3=\xi \lrcorner \;
[\;(\xi \lrcorner \; \varphi ) \wedge (\xi \lrcorner  \;\varphi )\wedge  \varphi \;]=\xi \lrcorner \; (6 |\xi |^2 \mu )= 6 \mu_{\xi}
\end{equation}
where $\mu_{\xi}=\mu |_{V_{\xi}}$ is the induced orientation form on $V_{\xi}$.

\vspace{.05in}

{\Lem $J_{\xi}$ is compatible with $\omega_{\xi}$, and it is metric invariant.}

\proof Let  $u,v\in V_{\xi}$
\begin{eqnarray*}
\omega_{\xi}(J_{\xi}(u), v )&=&\omega_{\xi} ( u\times \xi ,v)=
\langle \psi(u \times \xi,v),\xi\rangle =\varphi (u \times \xi,v, \xi)  \;\;\;\;\;\;\;\;\;\mbox {by} \;(5)
\\
&=&-\varphi (\xi, \;\xi \times u,\;v)= -\langle \;\xi\times (\xi\times u), v \; \rangle \;\;\;\;\;\;\;\;\;\;\;\;\;\;\;\;\;\;\;\;\;\;\mbox {by} \;(3)\\
&=& -\langle \; - | \xi |^2 u + \langle \xi,u \rangle \xi, v \;\rangle =  | \xi |^2 \langle u,v\rangle - \langle \xi,u \rangle \langle \xi,v\rangle \;\;\;\;\;\mbox {by} \;(9)\\
&=& \langle u,v\rangle .
\end{eqnarray*}
By plugging in $J_{\xi}(u)$, $J_{\xi}(v)$ for $u$, $v$:
$\langle J_{\xi}(u), J_{\xi}(v)\rangle =-\omega_{\xi}(u,J_{\xi}(v))=\langle u,v\rangle$
\qed
\vspace{.05in}

{\Lem  $\Omega_{\xi} $ is a non-vanishing $(3,0)$ form.}

\proof  By a local calculation as in Section 3.1 we see that
$\Omega_{\xi} $ is a  $(3,0)$ form, and is non-vanishing bacause
$\Omega_{\xi} \wedge \overline{\Omega}_{\xi}=8i \;vol(X_{\xi})$,
i.e.
\begin{eqnarray*} \frac{1}{2i} \;\Omega_{\xi} \wedge \overline{\Omega}_{\xi}=Im\;\Omega_{\xi} \wedge Re\;\Omega_{\xi} &=& (\xi \lrcorner *\varphi)\wedge [\;\xi\lrcorner \;(\xi^{\#}\wedge \varphi)\;]\\
&=&-\xi \lrcorner \; [\;(\xi \lrcorner  *\varphi)\wedge (\xi^{\#}\wedge \varphi)\;]\\
&=&\xi \lrcorner \;[ *(\xi^{\#}\wedge \varphi )\wedge ( \xi^{\#}\wedge \varphi ) \;]\;\;\;\;\;\;\;\mbox {by} \;(10)\\
&=& |\xi^{\#}\wedge \varphi |^2 \; \xi \lrcorner \; vol (M)\\
&=& 4|\xi^{\#}|^2 \; (*\xi^{\#}) =4 \;vol (X_{\xi} ). \;\;\;\;\;\;\;\;\mbox {by} \;(7)
\end{eqnarray*} \qed

\vspace{.1in}

We can easily calculate  $*Re \;\Omega_{\xi}=-Im \;
\Omega_{\xi}\wedge \xi^{\#}$ and  $*Im \;\Omega_{\xi}=Re \;
\Omega_{\xi}\wedge \xi^{\#}$. In particular if  $\star $ is the
star operator of $X_{\xi}$ (for example by (15) $\star
\omega_{\xi}=\omega_{\xi}^2$/2), then
\begin{equation}
\star Re \; \Omega_{\xi}=Im\; \Omega_{\xi}.
\end{equation}

\vspace{.05in}

Notice that $\omega_{\xi} $ is a symplectic structure on $X_{\xi}$
whenever  $d \varphi =0$  and $\CL_{\xi}(\varphi)|_{V_{\xi}} =0$,
where $\CL_{\xi}$ denotes the Lie derivative along $\xi$. This is
because $\omega_{\xi} ={\xi} \lrcorner \;\varphi $ and:
$$ d\omega_{\xi}=\CL_{\xi}(\varphi) - \xi  \lrcorner \; d\varphi = \CL_{\xi}(\varphi). $$

\vspace{.1in}

\noindent Also $d^{*}\varphi =0 \implies d^{\star}\omega_{\xi}=0 $, without any condition on the vector field $\xi$,  since
\begin{equation}
*\varphi=\star \;\omega_{\xi}- Im \; \Omega_{\xi}\wedge \xi^{\#},
\end{equation}

\vspace{.05in}

\noindent and hence $d(\star \omega_{\xi})=d(*\varphi |_{X_{\xi}}
)=0$. Also $d\varphi=0 \implies $ $d
(Re\;\Omega_{\xi})=d(\varphi|_{X_{\xi}})=0$.

\vspace{.05in}

 Furthermore, $d^{*}\varphi =0$ and
$\CL_{\xi}(*\varphi)|_{V_{\xi}} =0 \implies $ $d( Im\;
\Omega_{\xi})=0 $; this is because $Im\;\Omega_{\xi} =
\xi\lrcorner \; (*\varphi)$, where $*$ is the star operator on
$(M, \varphi )$.  Also,  $J_{\xi}$ is integrable when $d\Omega=0$
(e.g. \cite{hi}). By using the following definition, we can sum up all the conclusions of the above discussion as Theorem 5 below.

\vspace{.05in}

{\Def $(X^6,  \omega, \Omega, J)$  is called an {\it almost
Calabi-Yau manifold}, if $X$ is a Riemannian manifold with a
non-degenerate $2$-form $\omega$ (i.e. $\omega^3 =6vol (X)$) which
is co-closed, and  $J$ is a metric invariant almost complex
structure which is  compatible with $\omega$, and $\Omega $ is a
non-vanishing $(3,0)$ form with  $Re\;\Omega $ closed.
Furthermore, when  $\omega$ and $ Im\; \Omega$ are closed,  we
call this a Calabi-Yau manifold.}


{\Thm Let $(M,\varphi)$ be a $G_2$ manifold, and $\xi $ be a  unit
vector field which comes  from a codimension one foliation on $M$,
then  $(X_{\xi},\omega_{\xi}, \Omega_{\xi},J_{\xi})$ is an almost
Calabi-Yau manifold with $\varphi |_{X_{\xi}}= Re\; \Omega_{\xi} $
and  $*\varphi |_{X_{\xi}}= \star \omega_{\xi} $. Furthermore, if
$\CL_{\xi}(\varphi )|_{X_{\xi}}=0$ then $d\omega_{\xi}=0$,  and if
$\CL_{\xi}(*\varphi)|_{X_{\xi}}=0$ then $J_{\xi}$ is integrable; when both of these conditions are satisfied then
$(X_{\xi},\omega_{\xi}, \Omega_{\xi},J_{\xi})$  is a Calabi-Yau
manifold.}

{\Rm If $\xi$ and $\xi'$ are sections of ${\bf V}$ and ${\bf E}$  respectively, then from \cite{m}  the condition
$\CL_{\xi}(*\varphi)|_{X_{\xi } }=0$  (complex geometry of $X_{\xi }$)
 implies that deforming associative submanifolds of $X_{\xi}$ along $\xi$  in $M$ keeps them associative; and
$\CL_{\xi '}(\varphi )|_{X_{\xi ' }}=0$ (symplectic geometry of $X_{\xi '}$) implies  that deforming coassociative submanifolds of $X_{\xi '}$ along $\xi '$ in $M$  keeps them coassociative (e.g. for an example see Example 1).}

\vspace{.1in}

Notice that both complex and symplectic structure of the
CY-manifold $X_{\xi}$ in Theorem 3 is determined by $\varphi$ when
they exist. Recall that (c.f. \cite{v})  elements $\Omega \in
H^{3,0}(X_{\xi},\C)$ along with topology of $X_{\xi}$ (i.e. the
intersection form of $H^3(X_{\xi},\Z)$) parametrize complex
structures on $X_{\xi}$ as follows: We compute the third betti
number $b_{3}(M)=2h^{2,1} +2 $ since
$$H^{3}(X_{\xi},\C)=H^{3,0}\oplus H^{2,1} \oplus H^{1,2} \oplus
  H^{0,3}=2(\C\oplus H^{2,1}). $$

Let $\{ A^{i}, B_{j} \}$ be a symplectic basis of $H_{3}(X,\Z)$,
$i=1,.., h^{2,1}+1$, then
 \begin{equation}
 X_{i}=\int_{A^{i}}\Omega
 \end{equation}
give complex numbers which are local homegenous coordinates of the
moduli space of complex structures on $X_{\xi}$, which is a
$h^{2,1}$ dimensional  space (there is an extra parameter here
since $\Omega$ is defined up to scale multiplication).
\vspace{.1in}

As we have seen in the  example of the Section 3.1, the choice of
$\xi$ can give rise to quite different complex structures on
$X_{\xi}$ (e.g. $SU(2)$ and $SU(3)$ structures).  For example,
assume $\xi \in \Omega^{0} (M, {\bf V})$ and $\xi'\in \Omega^{0}
(M, {\bf E})$ be unit vector fields, such that the the codimension
one plane fields  $\xi^{\perp}$ and $\xi'^{\perp}$ come from
foliations. Let $X_{\xi }$ and $X_{\xi '}$ be pages of the
corresponding foliations. By our definition $X_{\xi }$ and
$X_{\xi' }$ are mirror duals of each other. Decomposition
$T(M)={\bf E}\oplus{\bf V}$ gives rise to splittings
$TX_{\xi}={\bf E}\oplus \bar{ {\bf E} }$, and $TX_{\xi'}={\bf
C}\oplus  {\bf V} $, where $\bar{\bf E}= \xi^{\perp}( {\bf V}
)\subset {\bf V} $ is a $3$-dimensional subbundle, and ${\bf C}=(
\xi' )^{\perp} ({\bf E})\subset {\bf E} $ is a  $2$-dimensional
subbundle. Furthermore, ${\bf E}$ is Lagrangian in $TX_{\xi}$ i.e.
$J_{\xi}({\bf E})=\bar{\bf E}$, and ${\bf C}$, ${\bf V}$ are
complex in $TX_{\xi'}$  i.e. $J_{\xi'}({\bf C})={\bf C}$ and
$J_{\xi'}({\bf V})={\bf V}$.
Also notice that, $Re\;\Omega_{\xi}$ is a
calibration form of  ${\bf E}$ , and $\omega_{\xi}$ is a
calibration form of ${\bf C}$. In particular,  $\langle
\Omega_{\xi}, {\bf E})=1$ and $\langle \omega_{\xi }\wedge \xi^{
\#}, {\bf E}\rangle=0$; and $\langle \Omega_{\xi'}, {\bf E})=0$
and   $\langle \omega_{\xi' }\wedge (\xi')^{ \#}, {\bf
E}\rangle=1$.

\vspace{.1in}

If $X_{\xi}$ and $X_{\xi'}$ are strong duals of each other, we can
find a homotopy of non-vanishing unit vector fields $\xi_{t}$
($0\leq t \leq 1$) starting with  $\xi\in {\bf V}$ ending with
$\xi'\in {\bf E}$. This gives a $7$-plane distribution
$\Xi=\xi^{\perp}_{t} \oplus \frac{\partial}{\partial t}$ on
$M\times [0,1]$ with integral submanifolds $X_{\xi} \times
[0,\epsilon)$ and $X_{\xi'} \times (1-\epsilon, 1]$ on a
neighborhood of the boundary. Then by \cite{th} and \cite{th1} we
can homotop $\Xi$ to a foliation  extending the foliation on the
boundary (possibly by taking $\epsilon$ smaller). Let $Q^7\subset
M\times [0,1]$ be the smooth manifold given by this foliation,
with $\partial Q= X_{\xi} \cup X_{\xi'}$, where $X_{\xi}\subset
M\times \{0\}$ and $X_{\xi'}\subset M\times \{1\}$.

\vspace{.1in}

  \begin{figure}[ht]  \begin{center}
\includegraphics{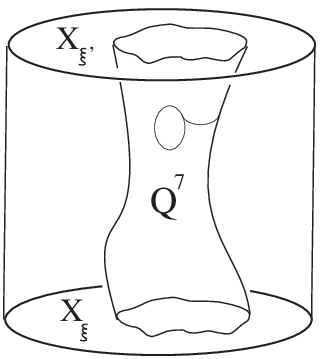}   \caption{}    \end{center}
   \end{figure}

\vspace{.1in}

 We can define $\Phi \in \Omega^{3}(M \times [0,1])$
with $\Phi |_{X_{\xi}}= \Omega_{\xi}$ and  $\Phi |_{X_{\xi'}} =\xi
\lrcorner \star \omega_{\xi'} $
 \begin{equation*}
\Phi =\Phi (\varphi, \Lambda, t )= \langle   \omega_{\xi_{t}}\wedge \xi_{t}^{\#}, {\bf E} \rangle \; \xi_{t}'' \lrcorner \star \omega_{\xi_{t}}
+ \langle Re\; \Omega_{\xi_{t}},  {\bf E}  \rangle \; \Omega_{\xi_{t}}
\end{equation*}
\noindent where $\xi''_{t}= J_{\xi\times \xi'}(\xi_{t})=\xi_{t}\times (\xi\times \xi' )$ (hence $\xi''_{0}=-\xi'$ and $\xi''_{1}=\xi$). This can be viewed as a correspondence between the complex structure of $X_{\xi}$ and the  symplectic structure of $X_{\xi'}$.
In general, the manifold pairs $X_{\alpha }$ and $X_{\beta }$  (as constructed in Theorem 5)  determine each others  almost Calabi -Yau structures via $\varphi$ provided they are defined.

\vspace{.1in}

{\Prop Let  $\{ \alpha, \beta \}$ be orthonormal vector fields
on $(M, \varphi) $. Then on $X_{\alpha}$ the following hold
\vspace{.05in}
\begin{itemize}
\item [(i)]$ Re\;\Omega_{\alpha }=\omega_{\beta} \wedge \beta^{\#} + Re\; \Omega_{\beta} $\\
\item [(ii)]$ Im\;\Omega_{\alpha }=\alpha  \lrcorner \; (\star \omega_{\beta})-(\alpha  \lrcorner \; Im\; \Omega_{\beta} )\wedge \beta^{\#} $\\
\item[(iii)] $\omega_{\alpha }= \alpha  \lrcorner \;Re \; \Omega_{\beta } + (\alpha  \lrcorner \; \omega_{\beta })\wedge \beta^{\#} $
\end{itemize} }

\proof  Since $Re\;\Omega_{\alpha } =\varphi |_{X_{\alpha}}$  (i) follows.
Since $Im\; \Omega_{\alpha}=\alpha  \lrcorner *\varphi$  following gives (ii)
\begin{eqnarray*}
\alpha  \lrcorner \; (\star \omega_{\beta})&=& \alpha  \lrcorner \; [\; \beta  \lrcorner \;*(\beta  \lrcorner \;\varphi )\;] \\
&=&\alpha  \lrcorner \; \beta  \lrcorner \; (\beta^{\#} \wedge * \varphi )\\
&=& \alpha  \lrcorner  *\varphi + \beta^{\#}\wedge (\alpha  \lrcorner \;\beta  \lrcorner  *\varphi) \\
&=&  \alpha  \lrcorner  *\varphi + (\alpha  \lrcorner \; Im\; \Omega_{\beta} )\wedge \beta^{\#}
\end{eqnarray*}

(iii) follows from the following computation
$$\alpha  \lrcorner \; Re\; \Omega_{\beta}= \alpha  \lrcorner \;\beta  \lrcorner \; (\beta ^{\#} \wedge \varphi)= \alpha  \lrcorner \; \varphi +\beta^{\#}\wedge (\alpha  \lrcorner \;\beta  \lrcorner \;\varphi )= \alpha  \lrcorner \; \varphi- (\alpha  \lrcorner \; \omega_{\beta })\wedge \beta^{\#} $$
\qed

\vspace{.1in}

Notice that even though the identities of Proposition 6 hold only
after restricting the right hand side to $X_{\alpha}$, all the
individual terms are defined everywhere on $(M,\varphi)$. Also,
from the construction, $X_{\alpha }$ and  $X_{\beta }$ inherit
vector fields $\beta$ and $\alpha$, respectively.

\vspace{.05in}

{\Cor Let  $\{ \alpha, \beta \}$ be orthonormal vector fields on
$(M,\varphi)$. Then there are  $A_{\alpha \beta} \in \Omega^{3}(M)$, and $W_{\alpha \beta}\in \Omega^{2}(M)$ satisfying
\begin{eqnarray*}
(a)\;\;\;\; \;\;\; \; \;\;\; \varphi |_{X_{\alpha}}= Re\;\Omega_{\alpha} &\mbox {and} & \varphi |_{X_{\beta}}=Re\; \Omega_{\beta}\\
(b)\;\;\;\; \;\; A_{\alpha \beta}|_{X_{\alpha}}= Im\;\Omega_{\alpha} &\mbox {and} &
A_{\alpha \beta}|_{X_{\beta}}= \alpha  \lrcorner \; (\star \omega_{\beta})\\
(c) \hspace{.5in} W_{\alpha \beta}|_{X_{\alpha}}= \omega_{\alpha} &\mbox {and} &
W_{\alpha \beta}|_{X_{\beta}}= \alpha  \lrcorner \;Re \; \Omega_{\beta }
\end{eqnarray*} }

For example, when $\varphi$ varies through metric preserving $G_2$ structures \cite{b2}, (hence fixing the orthogonal frame $\{\xi, \xi' \}$),
it induces variations of $\omega$ one side, and $\Omega$ on the
other side.

{\Rm By using Proposition 6,  in the previous torus  example of 3.1 one can show a natural correspondence between the groups  $H^{2,1}(X_{\xi} )$ and  $H^{1,1}(X_{\xi'} ) $.  Even though ${\mathbb T}^7$ is a trivial example of a $G_2$ manifold, it is an important special case since the $G_2$ manifolds of Joyce are obtained by smoothing quotients of ${\mathbb T}^7$ by finite group actions. We believe this process turns the subtori $X_{\xi}$'s into Borcea-Voisin manifolds with a similar correspondence of their cohomology groups. }

\vspace{.1in}

For the discussion of the previous paragraph to work, we need  a
non-vanishing vector field $\xi $ in $T(M) ={\bf E}\oplus {\bf
V}$, moving from ${\bf V}$ to ${\bf E}$. The bundle ${\bf E}$
always has a non-zero section, in fact it has a non-vanishing
orthonormal $3$-frame field;  but ${\bf V}$ may not have a
non-zero section. Nevertheless the bundle  ${\bf V}\to M$ does
have a non-vanishing section in the complement of a $3$-manifold
$Y\subset M $, which is a transverse self intersection of the zero
section. In \cite{as}, Seiberg-Witten equations of such
$3$-manifolds were related to associative deformations.  So we can
use these partial sections and, as a consequence  $X_{\xi}$ and
$X_{\xi '}$ may not be  closed manifolds. The following is a
useful example:

\vspace{.05in}

{\Ex  Let $X_1$, $X_2$ be two Calabi-Yau manifolds,
where $X_1$ is the cotangent bundle of $S^3$ and $X_2$ is the
$\mathcal{O}(-1)\oplus\mathcal{O}(-1)$ bundle of $S^2$. They are conjectured to be the mirror duals of each other by physicists (c.f. \cite{ma}). By using the approach of this paper, we identify them as 6-dimensional submanifolds of a $G_2$ manifold. Let's choose $M=\wedge ^2_+(S^4)$; this is a
 $G_2$ manifold by Bryant-Salamon \cite{bs}.

\vspace{.05in}

Let $\pi: \wedge ^2_+ (S^4) \to S^4 $ be the bundle projection. The sphere bundle of $\pi$ (which is also $\overline{{\C\P}^3} $) is the so-called {\it twistor bundle}, let us denote it by $\pi_{1}: Z(S^4) \to S^4$. It is known that the normal bundle of each fiber $\pi_{1}^{-1}(p)\cong S^2$ in $Z(S^4)$ can be identified by
 $\mathcal{O}(-1)\oplus\mathcal{O}(-1)$ \cite{s}. Now we take ${\bf E}$ to be the bundle of vertical tangent vectors of $\pi$, and
${\bf V}=\pi^{*} (T S^4)$, lifted by connection distribution. Let $\xi$ be the pull-back of the vector field on $S^4$ with two zeros (flowing from  north pole ${\bf n}$  to south pole ${\bf s}$),  and let $\xi '$ be the radial vector field of ${\bf E}$ . Clearly $X_{\xi }=T^{*}(S^3)$ and $X_{\xi' }=\mathcal{O}(-1)\oplus\mathcal{O}(-1)$.
\vspace{.05in}

Note that  $\xi$ is non-vanishing in the complement of  $\pi^{-1}\{{\bf n,s}\}$, whereas  $\xi'$ is non-vanishing in the complement of the zero section of $\pi $. Clearly on the set where they are both defined,  $\xi$ and $\xi'$ are homotopic through nonvanishing vector fields $\xi_{t}$. This would define a cobordism between the complements of the zero sections of the bundles  $T^{*}(S^3)$ and $\mathcal{O}(-1)\oplus\mathcal{O}(-1)$, if the distributions $\xi_{t}^{\perp}$ were involutive.}

  \begin{figure}[ht]  \begin{center}
\includegraphics{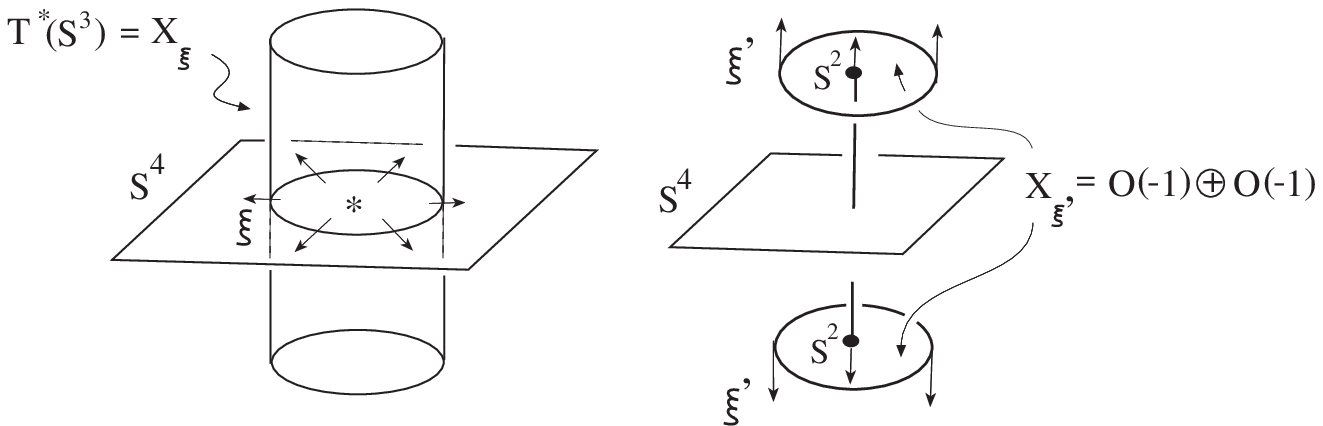}   \caption{}    \end{center}
   \end{figure}
\vspace{.1in}

\vspace{.1in}

Here the change of complex structures $X_{\xi '} \leadsto X_{\xi }$ happens as follows. Let $S^{3}_{\lambda}\to S^2$ be the Hopf map with fibers consisting of circles of radius
$\lambda $, clearly $S^{3}_{\infty}=S^2\times \R$
$$(\C^{2}- 0 )\times S^2 \to (S_{\lambda}^3 \times \R)\times S^{2}\; \stackrel {\lambda \to \infty}{\longrightarrow} \;
 S^3_{\infty}\times S^3_{\infty}$$

\noindent where the complex structure on $ S^3_{\infty}\times S^3_{\infty}$ is the obvious one, induced from exchanging the factors. In general if we allow the vector fields $\xi$ and $\xi '$ be homotopic through vector fields $\xi_{t}$ possibly with zeros, or the family $\xi_{t}^{\perp}$ not remain involutive  the cobordism between $X_{\xi}$ and $X_{\xi'}$ will have singularities.

\vspace{.1in}

{\Rm If we apply the construction of Example 1 to the total space of the spinor bundle $Q\to S^3$ (see Remark 2), the two dual $6$-manifolds we get are $S^2\times {\R}^4$ and $S^3\times \R^3$.}

\vspace{.1in}

There is also a concept of mirror-dual $G_2$ manifolds in a Spin(7) manifold,  hence we can talk about  mirror dual CY manifolds coming from two different  mirror dual $G_2$ submanifolds of a Spin(7) manifold. This is the the subject of the next section.

\vspace{.05in}

\section{ Mirror duality in $Spin(7)$ Manifolds}

\vspace{.1in}

Similar to Calabi-Yau case there is a notion of mirror duality
between $G_2$ manifolds \cite{ac}, \cite{av}, \cite{ gyz},
\cite{sv}. In this section we will give a definition of mirror
$G_2$ pairs, and an example which shows that associative and
co-associative geometries in mirror $G_2$ pairs are induced from
the same calibration $4$-form in a $Spin(7)$ manifold, and hence
these geometries dual to each other. Let us first recall the basic
definitions and properties of $Spin(7)$ geometries. The main
references in this subject are \cite{hl} and \cite{ti}.

{\Def An 8-dimensional Riemannian manifold $(N,\Psi)$ is called a
{\em Spin(7) manifold} if the holonomy group of its metric
connection lies in $Spin(7)\subset GL(8)$.}

\vspace{.05in}

Equivalently, a $ Spin(7) $ manifold is an $8$-dimensional
Riemannian manifold with a triple cross product $\times $ on its
tangent bundle, and a closed $4$-form $\Psi \in \Omega^{4}(N)$
with
$$ \Psi (u,v,w,z)=\langle u \times v \times w,z \rangle.$$
{\Def  A 4-dimensional submanifold $X$ of a $Spin(7)$ manifold $(N,\Psi)$ is
called {\em Cayley} if $\Psi|_X\equiv vol(X)$. }

\vspace{.05in}

Analogous to the $G_2$ case, we introduce a tangent bundle valued
$3$-form, which is just the triple cross product of $N$.

{\Def Let $(N, \Psi )$ be a $Spin(7)$ manifold. Then $\Upsilon \in
\Omega^{3}(N, TN)$ is the tangent bundle valued 3-form defined by
the identity:}
\begin{equation*}
\langle \Upsilon (u,v,w) , z \rangle=\Psi  (u,v,w,z)=\langle
u\times v\times w , z \rangle .
\end{equation*}

$Spin(7)$ manifolds can be constructed from $G_2$
manifolds. Let $(M,\varphi)$ be a $G_2$ manifold with a 3-form
$\varphi$, then $M\times S^1$ (or $M\times \mathbb{R}$) has holonomy group $G_2\subset Spin(7)$, hence is
a Spin(7) manifold. In this case $\Psi= \varphi\wedge dt +
*_{7} \varphi$, where $*_{7}$ is the star operator of $M^7$.

\vspace{.05in}

Now we will repeat a similar construction for  a $Spin(7)$ manifold $(N, \Psi )$, which we did for $G_2$ manifolds. Here we make an assumption that $T(M)$ admits a nonvanishing $3$-frame field $\Lambda =\langle u,v,w \rangle$, then we decompose $T(M)={\bf K}\oplus{\bf D}$, where ${\bf K}=\langle u,v,w, u\times v\times w \rangle$ is the bundle of Cayley $4$-planes  (where $\Psi$ restricts to be 1) and ${\bf D}$ is the complementary subbundle (note that this is also a bundle of Cayley $4$-planes since the form $\Psi$ is self dual).  In the $G_2$ case, existence of an analogous decomposition of the tangent bundle followed from \cite{t} (in this case we can just restrict to a submanifold which a $3$-frame field exists).
On a chart in $N$ let  $e_1,...e_8$ be an orthonormal frame and $e^1,...,e^8$ be the dual coframe, then  the calibration
4-form is given as (c.f. \cite{hl}) \:
\begin{equation}
\begin{aligned}
\Psi=\;\; e^{1234}&+(e^{12}-e^{34})\wedge (e^{56}-e^{78})\\
&+(e^{13}+e^{24})\wedge (e^{57}+e^{68})\\
&+(e^{14}-e^{23})\wedge (e^{58}-e^{67})+e^{5678}\\
\end{aligned}
\end{equation}

\noindent which is a self dual $4$-form, and the corresponding tangent bundle valued 3-form is

\vspace{.1in}

\begin{equation*}
\begin{aligned}
\Upsilon=&\;\;(e^{234}+e^{256}-e^{278}+e^{357}+e^{368}+e^{458}-e^{467})e_1\\
&+(-e^{134}-e^{156}+e^{178}+e^{457}+e^{468}-e^{358}+e^{367})e_2\\
&+(e^{124}-e^{456}+e^{478}-e^{157}-e^{168}+e^{258}-e^{267})e_3\\
&+(-e^{123}+e^{356}-e^{378}-e^{257}-e^{268}-e^{158}+e^{167})e_4\\
&+(e^{126}-e^{346}+e^{137}+e^{247}+e^{148}-e^{238}+e^{678})e_5\\
&+(-e^{125}+e^{345}+e^{138}+e^{248}-e^{147}+e^{237}-e^{578})e_6\\
&+(-e^{128}+e^{348}-e^{135}-e^{245}+e^{146}-e^{236}+e^{568})e_7\\
&+(e^{127}-e^{347}-e^{136}-e^{246}-e^{145}+e^{235}-e^{567})e_8.\\
\end{aligned}
\end{equation*}

\vspace{.1in}

This time we show that the form $\Upsilon$ induce $G_2$ structures
on certain subbundles of $T(N)$. Let $\gamma$ be a nowhere
vanishing vector field of $N$. We define a $G_2$ structure
$\varphi_{\gamma}$ on the $7$-plane bundle
$V_{\gamma}:=\gamma^{\perp}$ by (where $*_{8}$ is the star operator on $N^8$)
\begin{equation}
\varphi_\gamma:=\langle \Upsilon, \gamma \rangle = \gamma \lrcorner \; \Psi= *_{8}(\Psi\wedge \gamma^{\#}).
\end{equation}

 Assuming that $V_{\gamma}$ comes from a foliation, we let $M_{\gamma}$ be an integral submanifold of $V_{\gamma}$. We have $d\varphi_{\gamma}=0$, provided $\CL_{\gamma } (\Psi )|_{V_{\gamma}} =0 $. On the other hand, we always have  $d(\star\varphi_{\gamma})=0$ on $M_{\gamma}$. To see this, we use
 $$ \Psi =\varphi_{\gamma}\wedge \gamma^{\#} +*_{7} \varphi_{\gamma}$$ where $*_{7} $ is the star operator on $M{\gamma}$,  and use $d\Psi=0$ and  the foliation condition $d \gamma^{\#}\wedge \gamma^{\#}=0$, and the identity $\theta|_{M_{\gamma}}= \gamma \lrcorner \; [\;\theta \wedge \gamma^{\#}\;]$ for  forms $\theta $.  In order to state the next theorem, we need a  definition:

 \vspace{.05in}
 {\Def A manifold with  $G_2$ structure $(M,\varphi )$ is called an almost $G_2$-manifold if $\varphi $ is co-closed.}

\vspace{.05in}
{\Thm Let $(N^8, \Psi )$ be a $Spin(7)$ manifold, and $\gamma $ be a unit vector field which comes from a foliation, then $(M_{\gamma},\varphi_{\gamma})$ is an almost $G_2$ manifold. Furthermore if $\CL_{\gamma} (\Psi)|_{M_{\gamma}}=0$ then $(M_{\gamma},\varphi_{\gamma})$  is a $G_2$ manifold.}
\proof Follows by checking Definition 1 and by the discussion above.\qed

\vspace{.1in}

The following theorem says the induced $G_2$ structures on
$M_{\alpha}$, $M_{\beta }$ determine each other via $\Psi$; more
specifically $\varphi_{\alpha}$ and $\varphi_{\beta}$  are
restrictions of a global $3$-form of $N$.

\vspace{.05in}

{\Prop Let $(N, \Psi )$ be a $Spin(7)$ manifold, and $\{ \alpha,
\beta \}$ be an orthonormal vector fields on $N$.Then the
following holds on $M_{\alpha}$
$$ \varphi_{\alpha} = - \alpha \lrcorner  \; (\varphi_{\beta} \wedge \beta^{\#} +*_{7}\;\varphi_{\beta})$$}

\proof The proof follows from the definitions, and by expressing
$\varphi_{\alpha}$ and $\varphi_{\beta}$ in terms of $\beta^{\#}$
and $\alpha^{\#}$ by the formula (13). \qed

\vspace{.1in}

As in the $G_2$ case,  by choosing  different $\gamma $'s, one can
find various different $G_2$ manifolds  $M_{\gamma}$ with
interesting structures. Most interestingly, we will get certain
``dual'' $M_{\gamma} $'s by choosing $\gamma $ in ${\bf K}$ or in
${\bf D}$. This will shed light on more general version of mirror
symmetry of Calabi-Yau manifolds. First we will discuss an
example.

\vspace{.1in}

\subsection {An example}

$\;$
\vspace{.1in}

Let  $\mathbb
{T}^8=\mathbb{T}^4\times \mathbb{T}^4$ be the Spin(7) 8-torus, where
$\{e_1,e_2,e_3,e_4\}$ is the basis for the Cayley $\mathbb{T}^4$
and $\{e_5,e_6,e_7,e_8\} $ is the basis for the complementary
$\mathbb{T}^4$ . We can take
the corresponding calibration 4-form (20) above,
and take the decomposition $T(N)={\bf K}\oplus {\bf D}$,
where $\{e_1,e_2,e_3,e_4\}$ is the orthonormal basis for the
Cayley bundle ${\bf K}$, and $\{e_5,e_6,e_7,e_8\}$ is the local
orthonormal basis for the complementary bundle ${\bf D}$. Then if
we choose $\gamma=e_4=e_1\times e_2\times e_3$ then we get
$$\varphi_{\gamma}=
-e^{123}+e^{356}-e^{378}-e^{257}-e^{268}-e^{158}+e^{167}$$ On the
other hand, if we choose $\gamma'=e_5$ then we get
$$\varphi_{\gamma '}=e^{126}-e^{346}+e^{137}+e^{247}+e^{148}-e^{238}+e^{678}$$
which give different $G_2$ structures on the $7$ toris $M_{\gamma}$ and $M_{\gamma'}$.

\vspace{.1in}

Note that if we choose $\gamma$ from the Cayley bundle ${\bf K}$,
we get the $G_2$ structure on the 7-torus $M_{\gamma}$ which reduces the Cayley
4-torus $\mathbb{T}^4=\mathbb{T}^3 \times S^{1} $ (where $\gamma $ is tangent to  $S^1$ direction) to an associative 3-torus $ \mathbb{T}^3 \subset M_{\gamma}$ with respect to this $G_2$
structure. On the other hand if we choose $\gamma '$ from the
complementary bundle ${\bf D}$, then the Cayley 4-torus $\mathbb{T}^4$ will be a coassociative submanifold of the $7$-torus $M_{\gamma'}$with the corresponding $G_2$ structure. Hence associative and coassociative geometries are dual to each other as they are induced from the same
calibration 4-form $\Psi$ on a Spin(7) manifold. This suggests the
following definition of the ``mirror duality'' for $G_2$
manifolds.

\vspace{.05in}

{\Def Two $7$-manifolds with $G_{2}$ structures are mirror pairs, if their $G_2$-structures are induced from the same calibration 4-form in a
$Spin(7)$ manifold. Furthermore they are strong duals if their normal vector fields  are homotopic. }

\vspace{.05in}

{\Rm For example, by \cite{bs} the total space of an $\R^4 $ bundle over
$S^4$ has a $Spin(7)$ structure. By applying the process of Example 1, we obtain mirror pairs $M_{\gamma}$ and $M_{\gamma'}$ to be $S^3\times \R^4$ and  $\R^4\times S^3$ with dual  $G_2$ structures.}

\vspace{.1in}

\subsection{Dual Calabi-Yau's inside of Spin(7)}

$\;$
\vspace{.1in}

Let $(N^8, \Psi)$ be a $Spin(7)$ manifold, and let $\{\alpha, \beta\}$ be an orthonormal  $2$-frame field in $N$, each coming from a foliation
 Let $(M_{\alpha}, \varphi_{\alpha})$ and $(M_{\beta}, \varphi_{\beta})$ be the $G_2$ manifolds given by Theorem 8. Similarly by Theorem 5, the vector fields  $\beta $ in $M_{\alpha }$, and $\alpha $ in $M_{\beta }$ give   almost Calabi-Yau's $X_{\alpha \beta}\subset M_{\alpha}$ and $X_{\beta \alpha }\subset M_{\beta }$. Let us denote $X_{\alpha \beta}= (X_{\alpha\beta}, \omega_{\alpha \beta}, \Omega_{\alpha  \beta}, J_{\alpha \beta})$ likewise
$X_{\beta \alpha }= (X_{\beta \alpha}, \omega_{\beta \alpha }, \Omega_{ \beta \alpha }, J_{\beta, \alpha})$. Then we have

\vspace{.1in}

{\Prop The following relations hold: \vspace{.05in}
\begin{itemize}
\item [(i)] $J_{\alpha \beta}(u)=u\times \beta\times \alpha $\\
\item [(ii)] $\omega_{\alpha \beta}=\beta \lrcorner \; \alpha \lrcorner \; \Psi$\\
\item [(iii)] $Re\; \Omega_{\alpha \beta}=\alpha \lrcorner \;
\Psi |_{X_{\alpha \beta} }$\\
\item [(iv)] $Im\; \Omega_{\alpha \beta}=\beta \lrcorner \;
\Psi |_{X_{\alpha \beta} }$
\end{itemize}}
\proof (i), (ii), and (iii) follow from definitions, and the formula $X\times Y= (X \lrcorner  \; Y\lrcorner  \;\varphi)^{\#}$.
\begin{eqnarray*}
Im\; \Omega_{\alpha \beta}&=& \beta\lrcorner  *_{7} \varphi_{\alpha}=
 \beta\lrcorner  *_{7} (\alpha \lrcorner \; \Psi)\\
 &=&\beta \lrcorner  \;[\alpha \lrcorner  \;*_{8}(\alpha \lrcorner  \; \Psi)]\\
 &=&\beta \lrcorner  \; [\alpha \lrcorner  \; (\alpha^{\#}\wedge \Psi )]\;\;\;\;\;\;\;\;\;\;\mbox{by}\;(10)\\
 &=&\beta \lrcorner  \; [\Psi- \alpha^{\#}\wedge (\alpha  \lrcorner  \; \Psi)] \\
  &=&\beta  \lrcorner  \; \Psi -\alpha^{\#}\wedge (\beta  \lrcorner  \; \alpha  \lrcorner  \; \psi)] .
\end{eqnarray*}

Left hand side is already defined on  $X_{\alpha \beta}$, by restricting to $X_{\alpha \beta}$ we get (iii). \qed

\vspace{.05in}

{\Cor When $X_{\alpha, \beta}$ and  $X_{\beta, \alpha}$ coincide, they are oppositely oriented manifolds and $\omega_{\alpha, \beta}=- \omega_{\beta, \alpha}$, and $Re\; \Omega_{\alpha \beta}=-Im\; \Omega_{\beta \alpha}$ (as forms on $X_{\alpha \beta}$).}

\vspace{.1in}

Now let   $\{\alpha, \beta, \gamma\}$ be  an orthonormal $3$-frame field in $(N^8, \Psi)$, and let
$(M_{\alpha}, \varphi_{\alpha})$, $(M_{\beta}, \varphi_{\beta})$,  and $( M_{\gamma}, \varphi_{\gamma})$ be the corresponding almost $G_2$ manifolds. As before, the orthonormal vector fields  $\{ \gamma, \beta \}$ in $M_{\alpha }$ and $\{\gamma, \alpha \}$ in $M_{\beta }$ give rise to corresponding almost Calabi-Yau's  $X_{\alpha,\gamma}$, $X_{\alpha,\beta }$ in $M_{\alpha}$, and $X_{\beta,\gamma }$, $X_{\beta,\alpha }$ in $M_{\beta }$.

\vspace{.05in}

In this way $(N^8,\Psi )$ gives rise to $4$  Calabi-Yau descendents. By Corollary 11, $X_{\alpha \beta}$ and $X_{\beta \alpha}$ are different geometrically; they may not even be the same as smooth manifolds, but for simplicity we may consider it to be the same smooth manifold obtained from the triple intersection of the three $G_2$ manifolds.

\vspace{.05in}

In case we have a decomposition $T(N)={\bf K}\oplus {\bf D}$ of the tangent bundle of $(N^8,\Psi )$ by Cayley plus its orthogonal bundles(Section 4);  we can choose our frame special and obtain  interesting CY-manifolds. For example,  if we choose  $\alpha \in \Omega^{0}(M,{\bf K})$ and $\beta, \gamma \in \Omega^{0}(N,{\bf D})$ we get one set of complex structures, whose types are  indicated by the first row of the following diagram. On the other hand, if we choose all $\{\alpha, \beta,\gamma \}$ lie entirely in ${\bf K}$ or ${\bf D }$ we get an other set of complex structures, as indicated by the second row of the diagram.

\vspace{.1in}

$\begin{array}{ccccccccc}
&&&&(N^8, \Psi )&&&&\\
&&&\swarrow &&\searrow&&&\\
&&(M_{\alpha},\varphi_{\alpha})&&&&(M_{\beta}, \varphi_{\beta})&&\\
&\swarrow&&\searrow&&\swarrow&&\searrow&\\
{\bf X}_{\alpha \gamma}&&&&X_{\alpha \beta}\;\;\;\;\;\;\;{\bf X}_{\beta \gamma} &&&&X_{\beta \alpha}\\
&&&&&&&&\\
SU(3)&&&&SU(3)\;\;\;\;\;\;\; SU(2) &&&& SU(3)\\
&&&&&&&&\\
SU(2)&&&&SU(2)\;\;\;\;\;\;\; SU(2) &&&& SU(2)\\
\end{array}$

\vspace{.1in}

Here all the corresponding symplectic and the holomorphic forms of the resulting Calabi-Yau's come from restriction of global forms induced by $\Psi $. The following gives  relations between the complex/symplectic structures of these induced CY-manifolds; i.e. the structures $X_{\alpha \gamma}$, $X_{\beta \gamma}$  and $X_{\alpha \beta }$ satisfy a certain triality relation.

\vspace{.1in}

{\Prop We have the following relations;

 \vspace{.05in}
\begin{itemize}
\item [(i)]$ Re\;\Omega_{\alpha \gamma}=\alpha \lrcorner (\; *_{6} \;\omega_{\beta \gamma }) +\omega_{\alpha \beta}\wedge\beta^{\#} $\\
\item [(ii)]$ Im\;\Omega_{\alpha \gamma}=\omega_{\beta\gamma }\wedge \beta^{\#} -\gamma \lrcorner \;*_{6}(\omega_{\alpha \beta} )$\\
\item[(iii)] $\omega_{\alpha \gamma }= \alpha  \lrcorner \; Im\;\Omega_{\beta \gamma }+ (\gamma \lrcorner \;\omega_{\alpha \beta} )\wedge \beta^{\#} $
\end{itemize} }

\vspace{.05in}

\noindent First we need to prove a lemma;

\vspace{.05in}

{\Lem The following relations hold;
\begin{eqnarray*}
\alpha  \lrcorner \;  *_{6}(\omega_{\beta \gamma})&=& \alpha  \lrcorner \;\Psi +\gamma^{\#}\wedge (\alpha  \lrcorner \;\gamma \lrcorner \; \Psi) + \beta^{\#}\wedge (\alpha  \lrcorner \;\beta \lrcorner \; \Psi)
 -\gamma^{\#}\wedge \beta^{\#}\wedge (\alpha  \lrcorner \; \gamma  \lrcorner \;  \beta  \lrcorner \;  \Psi). \\
  Im\;\Omega_{\beta \gamma}&=&-\gamma \lrcorner \;\Psi - \beta^{\#}\wedge (\gamma \lrcorner \;\beta \lrcorner \;\Psi).\\
Re\; \Omega_{\alpha \gamma }&=& \alpha  \lrcorner \;  \Psi- \gamma^{\#}\wedge (\gamma \lrcorner \; \alpha \lrcorner \;\Psi).
 \end{eqnarray*}
 \proof
  \begin{eqnarray*}
\alpha  \lrcorner \;  *_{6}(\omega_{\beta \gamma})  &=& \alpha \lrcorner \; [\gamma \lrcorner \; \beta \lrcorner \; *_{8}\;(\gamma \lrcorner \; \beta \lrcorner \Psi) ]\\
 &=& - \alpha \lrcorner \;\gamma \lrcorner \; \beta \lrcorner \;(\gamma^{\#}\wedge \beta^{\#}\wedge \Psi)  \\
 &=&-\alpha \lrcorner \; \gamma \lrcorner \; [ -\gamma^{\#} \wedge \Psi + \gamma^{\#}\wedge \beta^{\#}\wedge (\beta \lrcorner \; \Psi) ]  \\
  &=& \alpha  \lrcorner \;\Psi +\gamma^{\#}\wedge (\alpha  \lrcorner \;\gamma \lrcorner \; \Psi) + \beta^{\#}\wedge (\alpha  \lrcorner \;\beta \lrcorner \; \Psi) \\
 & & -\gamma^{\#}\wedge \beta^{\#}\wedge (\alpha  \lrcorner \; \gamma  \lrcorner \;  \beta  \lrcorner \;  \Psi) .
\end{eqnarray*}
\begin{eqnarray*}
 Im\;\Omega_{\beta \gamma}&=& \gamma  \lrcorner \; *_{7}\;(\beta  \lrcorner \; \Psi) \\
 &=& \gamma \lrcorner \; [\beta \lrcorner \; *_{8}\;(\beta  \lrcorner \; \Psi) ]\\
 &=& - \gamma \lrcorner \; \beta \lrcorner \; (\beta^{\#}\wedge \Psi) \\
 &=& -\gamma \lrcorner \; \Psi - \beta ^{\#}\wedge (\gamma \lrcorner \; \beta \lrcorner \; \Psi).
 \end{eqnarray*}
 $$Re\; \Omega_{\alpha \gamma }=(\alpha  \lrcorner \;  \Psi ) |_{X_{\alpha \gamma}}= \gamma  \lrcorner \; [ \gamma^{\#} \wedge (\alpha  \lrcorner \; \Psi) ]= \alpha  \lrcorner \;  \Psi- \gamma^{\#}\wedge (\gamma \lrcorner \; \alpha \lrcorner \;\Psi).$$
\qed

\noindent  {\it Proof of Proposition 12}.    We calculate the
following by using Lemma 13:
\begin{eqnarray*}
 \alpha  \lrcorner \;  *_{6}(\omega_{\beta \gamma})+\omega_{\alpha \beta}\wedge \beta^{\#}&=& \alpha  \lrcorner \;\Psi +\gamma^{\#}\wedge (\alpha  \lrcorner \;\gamma \lrcorner \; \Psi) + \beta^{\#}\wedge (\alpha  \lrcorner \;\beta \lrcorner \; \Psi) \\
 & & -\gamma^{\#}\wedge \beta^{\#}\wedge (\alpha  \lrcorner \; \gamma  \lrcorner \;  \beta  \lrcorner \;  \Psi) +(\beta  \lrcorner \;\alpha  \lrcorner \;\Psi)\wedge \beta^{\#}\\
  &= &\alpha  \lrcorner \;\Psi - \gamma^{\#}\wedge (\gamma  \lrcorner \;\alpha \lrcorner \; \Psi) -\gamma^{\#}\wedge \beta^{\#}\wedge (\alpha  \lrcorner \; \gamma  \lrcorner \;  \beta  \lrcorner \;  \Psi).
\end{eqnarray*}

\vspace{.05in}

Since we are restricting to $X_{\alpha \gamma}$ we can throw away
terms containing $\gamma^{\#}$ and get (i). We prove (ii)
similarly:

\vspace{.05in}

\begin{eqnarray*}
Im\;\Omega_{\alpha \gamma }&=&(\gamma \lrcorner \; *_{7}\varphi_{\alpha} )=
\gamma \lrcorner \; [\;\alpha \lrcorner \; *_{8}\varphi_{\alpha} \; ]\\
&=& \gamma \lrcorner \; [\alpha \lrcorner \; *_{8} (\alpha \lrcorner \;\Psi )]= - \gamma \lrcorner \; [\alpha \lrcorner \;  (\alpha^{\#} \wedge \;\Psi ) ) \\
&=&- \gamma \lrcorner \;\Psi  +
\gamma \lrcorner \;[ \alpha^{\#}\wedge (\alpha \lrcorner \; \Psi) ]\\
&=& - \gamma \lrcorner \;\Psi  -\alpha^{\#}\wedge( \gamma \lrcorner \;\alpha \lrcorner \; \Psi] .
\end{eqnarray*}
\begin{eqnarray*}
\omega_{\beta \gamma }\wedge \beta^{\#} - \gamma \lrcorner \;(*_{6}\;\omega_{\alpha \beta })&=& - \gamma \lrcorner \;(*_{6}\;\omega_{\alpha \beta })+ (\gamma \lrcorner \; \beta \lrcorner \; \Psi)\wedge \beta^{\#} \\
&=&-\gamma \lrcorner \;\Psi -\beta^{\#}\wedge (\gamma  \lrcorner \;\beta \lrcorner \; \Psi) -\alpha^{\#}\wedge (\gamma  \lrcorner \;\alpha \lrcorner \; \Psi) \\
 & & +\beta^{\#}\wedge \alpha^{\#}\wedge (\gamma  \lrcorner \; \beta  \lrcorner \;  \alpha  \lrcorner \;  \Psi)  + (\gamma \lrcorner \; \beta \lrcorner \; \Psi)\wedge \beta^{\#}.
\end{eqnarray*}

\vspace{.1in}

\noindent Here, we used Lemma 13 with different indices
$(\alpha,\beta, \gamma) \mapsto (\gamma, \alpha, \beta)$, and
since we are restricting to $X_{\alpha \gamma}$ we threw away
terms containing $\alpha^{\#}$. Finally, (iii) follows by plugging
in Lemma 13 to definitions. \qed

\vspace{.2in}

Following says that Calabi -Yau structures of $X_{\alpha \gamma}$ and $X_{\beta \gamma}$  determine each other via $\Psi$.  Proposition 14 is basically a consequence  of  Proposition 6 and Corollary 11.

\vspace{.1in}

{\Prop We have the following relations
\vspace{.05in}
\begin{itemize}
\item [(i)]$ Re\;\Omega_{\alpha \gamma}=\alpha \lrcorner \;( *_{6}\;\omega_{\beta \gamma }) - (\alpha  \lrcorner \; Re\;\Omega_{\beta \gamma}) \wedge \beta ^{\#}$\\
\item [(ii)]$ Im\;\Omega_{\alpha \gamma}=\omega_{\beta \gamma} \wedge \beta^{\#} +Im\; \Omega_{\beta \gamma}$\\
\item[(iii)] $\omega_{\alpha \gamma }= \alpha  \lrcorner \;Im \; \Omega_{\beta \gamma} + (\alpha  \lrcorner \; \omega_{\beta \gamma})\wedge \beta^{\#} $
\end{itemize} }

\proof All follow  from the definitions and Lemma 11 (and by ignoring $\alpha^{\#}$ terms). \qed


{\Rm After this paper was written we learned that the results similar to Theorem 5 already appeared in \cite{aw}, \cite{asa}, \cite{c}, \cite{cs} (where they use more restricted vector fields $\xi$), and we also  found out that the idea of studying induced hypersurface structures, from manifolds with exceptional holonomy  goes back to earlier works of Calabi and Gray. }


\vspace{.1in}


 \begin{thebibliography}{9999}


\bibitem[Ac]{ac} B.~Acharya, {\em On Mirror Symmetry for Manifolds of Exceptional Holonomy}, Nucl.Phys. B524 (1998) 269--282, hep-th/9707186.

\bibitem[AV]{av} M.~Aganagic and C.~Vafa, {\em $G_2$ Manifolds, Mirror Symmetry and Geometric Engineering}, hep-th/0110171.

\bibitem[AS]{as} S.~Akbulut and S.~Salur, {\em Deformations in $G_2$ manifolds}, math.GT/ 0701790

\bibitem[ASa]{asa} V.~Apostolov and S.~Salamon, {\em Kahler Reduction of metrics with holonomy $G_2$}, math.DG/0303197.

\bibitem[AW]{aw} M.~Atiyah and E.~Witten, {\em M-theory dynamics on a manifold of $G_2$ holonomy}, Adv. Theor. Math. Phys. 6 (2003),
1--106.

\bibitem[B1]{b1} R.L.~Bryant,  {\em Metrics with exceptional holonomy}, Ann. of Math 126 (1987), 525--576.

\bibitem[B2]{b2} R.L.~Bryant,  {\em Some remarks on $G_2$-structures}, GGT (2006),
math.DG/0305124 , v3.

\bibitem[BS]{bs} R.L.~Bryant and M.S.~Salamon, {\em On the construction of some complete metrics with exceptional holonomy},
Duke Math. Jour.,  vol. 58, no 3  (1989), 829--850.

\bibitem[C]{c} F.M.~Caberra, {\em SU(3)-structures on hypersurfaces of manifolds with $G_2$-structures},
math.DG/0410610.

\bibitem[CS]{cs} S.~Chiossi and S.~Salamon, {\em The intrinsic torsion of $SU(3)$ and $G_2$-structures}, math.DG/0202282.



\bibitem[GYZ]{gyz} S.~Gukov, S.T.~Yau and E.~Zaslow, {\em Duality and Fibrations on
$G_2$ Manifolds}, Turkish Jour. of Math., 27, (2003), 61--97.

\bibitem[Hi]{hi} N.J.~Hitchin,  {\em The moduli space of special Lagrangian submanifolds }, dg-ga/9711002.



\bibitem[HL]{hl} F.R. ~Harvey, and H.B. ~Lawson, {\em Calibrated geometries}, Acta. Math. 148 (1982), 47--157.


\bibitem[K]{k} S.~Karigiannis,  {\em Deformations of $G_2$ and $Spin(7)$ structures on manifolds}, math.DG/0301218.

\bibitem[J]{j} D.D.~Joyce, {\em Compact Manifolds with special holonomy}, OUP, Oxford, 2000.

\bibitem[M]{m} R.C.~McLean, {\em Deformations of calibrated submanifolds}, Comm. Anal. Geom.  6 (1998), 705--747.

\bibitem[Ma]{ma} M.~Marino, {\em Enumerative geometry and knot invariants}, Infinite dimensional groups and manifolds,  27--92, IRMA Lect. Math. Theor. Phys., 5, de Gruyter, Berlin, 2004. hep-th/0210145.

\bibitem[S]{s} S.~Salamon, {\em  Riemannian geometry and holonomy groups}, Pitman Res.Notes in Math. Series, no. 201.

\bibitem[SV]{sv} S.L.~Shatashvili and C.~Vafa {\em Superstrings and Manifolds of Exceptional Holonomy}, Selecta Math. 1 (1995) 347, hep-th/9407025.

\bibitem[T]{t} E.~Thomas, {\em Postnikov invariants and higher order cohomology operations}, Ann. of Math. vol 85 (1967), 184--217.

\bibitem[Th]{th} W.~Thurston, {\em Existence of codimension-one foliations}, Ann. of Math. (2) 104, (1976) 249--268

\bibitem[Th1]{th1} W.~Thurston, {\em Private communication}.

\bibitem[Ti] {ti} G.~Tian, {\it Gauge Theory and Calibrated Geometry I},
Ann. of Math. 151, (2000), 193--268.


\bibitem[V]{v} C.~Vafa, {\em Simmons Lectures (Harvard)}.

\end{thebibliography}
\end{document}